\documentclass[12pt]{amsart}
\usepackage{amsmath}
\usepackage{amssymb}
\usepackage{amscd}
\usepackage{amsthm}

\theoremstyle{plain}

\newtheorem{thm}[equation]{Theorem}

\newtheorem{prop}[equation]{Proposition}

\newtheorem{lem}[equation]{Lemma}
\newtheorem{cor}[equation]{Corollary}

\theoremstyle{remark}

\newtheorem{rem}[equation]{Remark}
\theoremstyle{definition}

\newtheorem{question}[equation]{Question}

\def\C{\mathbb C}
\def\R{\mathbb R}

\def\Z{\mathbb Z}

\def\N{\mathbb N}

\def\O{{\mathcal O}}
\newcommand{\acts}{\curvearrowright}

\newcommand{\al}{\alpha}

\def\eps{\epsilon}
\def\ga{\gamma}
\def\Ga{\Gamma}

\def\La{\Lambda}

\def\<{\langle}
\def\>{\rangle}

\def\be{\beta}

\newcommand{\ul}{\underline}

\begin{document}
\title{Arithmetic aspects of self-similar groups}
\author{Michael Kapovich}
\address{Department of Mathematics, University
of California, Davis, CA 95616, USA. kapovich@math.ucdavis.edu}
\thanks{This research was partially supported by the NSF grant
DMS-05-54349.}
\date{September 4, 2008}
 \maketitle

Let $\Ga$ be a group. A {\em virtual endomorphism} of $\Ga$ is a
homomorphism $\varphi: \La\to \Ga$, where $\La$ is a finite index
subgroup in $\Ga$. An {\em invariant normal subgroup} of $\varphi$
is a normal subgroup $N<\Ga$ contained in $\La$, so that
$\varphi(N)<N$. We will be interested in the question

\begin{question}
When does $\varphi$ have a nontrivial {\em invariant normal
subgroup}?
\end{question}

\medskip
We will only consider the case of lattices $\Ga$ in algebraic
(linear) Lie groups $G$. Therefore, without loss of generality,
the subgroups $\La$ can be taken torsion-free. In particular,
existence of a nontrivial invariant normal subgroup is equivalent
to the existence of an infinite invariant normal subgroups. Our
main result is

\begin{thm}\label{main}
Let $\Ga$ be an irreducible lattice  in a semisimple algebraic Lie
group $G$. Then the following are equivalent:

a. $\Ga$ is virtually isomorphic to an arithmetic lattice in $G$,
i.e., contains a finite index subgroup isomorphic to such
arithmetic lattice.

b. $\Ga$ admits a virtual endomorphism $\varphi: \La\to \Ga$
without nontrivial invariant normal subgroups.
\end{thm}

\medskip
We note that conjugacy classes of faithful self-similar actions of
a group $\Ga$ on a rooted tree (of finite valency) are in 1-1
correspondence with conjugacy classes  of virtual endomorphisms of
$\Ga$ which contain no nontrivial invariant normal subgroups
\cite{Nekrashevych,Nekrashevych-Sidky}. Therefore,

\begin{cor}
Let $\Ga$ be an irreducible lattice in a semisimple algebraic Lie
group $G$. Then the following are equivalent:

a. $\Ga$ is virtually isomorphic to an arithmetic lattice in $G$.

b. $\Ga$ admits an irreducible faithful self-similar action on a
regular  rooted tree (of finite valency).
\end{cor}

The only nontrivial ingredient in the proof of Theorem \ref{main}
is the following theorem of Margulis \cite{Margulis}:

\begin{thm}
An irreducible lattice $\Ga$ in an algebraic semisimple Lie group
$G$ is arithmetic if and only if its commensurator $Comm(\Ga)$ in
$G$ is non-discrete.
\end{thm}

\medskip
We start our proof of Theorem \ref{main} with some generalities
about arithmetic groups.

Let $F$ be an (algebraic) number field. Let $F_\infty$ denote the
completion of $F$ with respect to the valuation given by the
embedding $F\to \C$. Let $\O$ be the ring of integers of $F$.
Given a prime ideal $p$ in $\O$ we let $F_p$ denote the completion
of $F$ with respect to the nonarchimedean valuation on $F$
determined by $p$.

Let $\ul{G}$ be a semisimple algebraic group defined over $\O$.
 Then the group $G=\ul{G}(F_\infty)$ is a real
or complex semisimple Lie group. We will regard $G$ as the
isometry group of a symmetric space $X$ (the action has finite
kernel coming from the finite center of $G$). Similarly, the group
$G_p:= \ul{G}(F_p)$ acts as a group of isometries of a Euclidean
building $X_p$. The group of integer points $\ul{G}(\O_p)$ is the
stabilizer in $G_p$ of a special vertex $o_p\in X$.

Consider an irreducible arithmetic group $\Gamma\subset
\ul{G}(\O)$ in $G$. Recall \cite{Margulis} that the commensurator
$Comm(\Ga)$ of $\Ga$ is $\ul{G}(F)$, the group of $F$-points in
$\ul{G}$.

Given an element $\al\in Comm(\Ga)$ we let $\varphi=\varphi_\al$
denote the automorphism of $G$ induced by conjugation
$$
\varphi(g)= \al g \al^{-1}.
$$
Thus, $\al$ induces {\em virtual endomorphisms} of $\Ga$ defined
by taking sufficiently deep finite index %congruence
subgroups $\La\subset \Ga$, so that $\varphi(\La)\subset \Ga$, and
then restricting $\varphi$ to $\La$. (For instance, sufficiently
deep congruence subgroups will work.)

Since $\ul{G}(F)$ is dense in $G_p$, we can choose a prime ideal
$p$ and $\al\in Comm(\Ga)$ so that the corresponding isometry
$\al: X_p\to X_p$ has unbounded orbits (i.e. is hyperbolic).

\begin{rem}
Moreover, if $\al\in \ul{G}(F)$ induces an elliptic automorphism
of every Euclidean building $X_p$, then it belongs to a finite
extension of $\Ga$. However, we do not need this fact.
\end{rem}

\begin{lem}
Given $\al$ as above, for every finite index subgroup $\Ga'\subset
\Ga$, we have $\varphi(\Ga')\ne \Ga'$.
\end{lem}
\proof If $\al \Ga' \al^{-1}=\Ga'$, then $\Ga'$ and $\al$ generate
a discrete subgroup $\hat\Ga'\subset G$, which (since $\Ga'$ is a
lattice) necessarily has $\Ga$ as a finite index subgroup.
Therefore, the group $\hat\Ga$ generated by $\Ga$ and $\hat\Ga'$
is also a finite extension of $\Ga$. Since $\Ga$ fixes the point
$o_p\in X_p$, it follows that orbit $\hat\Ga\cdot o_p\subset X_p$
is bounded. This contradicts the fact that $\ga: X_p\to X_p$ is
hyperbolic. \qed

\begin{thm}\label{T1}
Given $\Ga, \La, \al$ as above, if $N<\La$ is a normal subgroup of
$\Ga$ so that $\varphi(N)\subset N$, then $N$ is finite.
\end{thm}
\proof We let $M:= X/N$. Then for every $N$ as above, we obtain an
isometric  covering
$$
q=\bar\al: M\to M
$$
induced by the endomorphism $\varphi$ of the fundamental group of
$M$. Note that the group $\Ga$ acts on $M$ isometrically, we use
the notation $\bar\ga$ for an isometry of $M$ induced by $\ga\in
\Ga$.

1. First, consider the case when $G$ has rank $\ge 2$. Then, by a
theorem of Margulis \cite{Margulis}, every normal subgroup of
$\Ga$ is either finite or has finite index in $\Ga$. Assuming $N$
has finite index in $\Ga$, the manifold $M$ has finite volume and,
since $q$ preserves the volume form, it has to be a
diffeomorphism. Thus, $\varphi(N)=N$, which contradicts our choice
of $\al$.
%Since $G$
%has trivial center, the group $\Ga$ cannot contain a nontrivial
%finite  normal subgroup.

\medskip
2. We now consider the more interesting case when $G$ has rank 1.
We assume that $N$ is an infinite group. Therefore, it is Zariski
dense in $G$.

Consider the iterations $q^k$ of the isometric endomorphism $q:
M\to M$. Let $M_{thick}, M_{thin}$ denote the thick and thin parts
of $M$ with respect to the Margulis constant of $X$. Pick a
connected  compact subset $C\subset M_{thick}$ whose fundamental
group maps onto a Zariski dense subgroup $H$ of $N$ and which
contains a fundamental domain for the action $\Ga\acts M_{thick}$.

 Then, by the Kazhdan-Margulis lemma, for every $k\ge 1$,
$q^k(C)$ is never contained in $M_{thin}$. Therefore, for each $k$
there exists $x_k\in C$ so that $q^k(x_k)\in M_{thick}$. Moreover,
since $q$ is isometric, there exists $\eps>0$ so that for every
$k\in \N$, the injectivity radius of $M$ on $q^k(C)$ is bounded
from below by $\eps$.

Since the action $\Ga\acts M_{thick}$ is cocompact (as $X/\Ga$ has
finite volume), for every $k$ there exists $\ga_k\in \Ga$ so that
$$
\bar\ga_k \circ q^k (x_k)\in C.
$$
Therefore, the sequence of isometries $\bar\be_k:= \bar\ga_k \circ
q^k$ is precompact.

\begin{lem}
The set $\bar{I}:= \{\bar\be_k, k\in \N\}$ is finite.
\end{lem}
\proof If not, then there will be arbitrarily large $k, m$ so that
$\bar\be_k\ne \bar\be_m$ and the restrictions $\bar\be_k|C,
\bar\be_m|C$ are arbitrarily close in the sup-metric. In
particular, for large $k, m$, they are homotopic and, hence,
induce the same (up to conjugation in $N$) map $H\to N$ given by
$$
h \mapsto \be_k h \be_k^{-1}, \quad h \mapsto \be_m h \be_m^{-1}.
$$

Since $H$ is Zariski dense in $G$, its centralizer in $G/Z(G)$ is
trivial and, thus we have the equality of the cosets
$$
\be_k\cdot N\cdot Z(G)=\be_m\cdot N\cdot Z(G).
$$
Since the center of $G$ acts trivially on $X$, we obtain that
$\bar\be_k=\bar\be_m$. Contradiction. \qed

\medskip
We continue with the proof of Theorem \ref{T1}. Let $I\subset
\ul{G}(F)$ denote the finite set of representatives of lifts of
the isometries $\bar\be_k\in \bar{I}$. Hence, for each $k\in \N$,
$$
\ga_k \cdot \al^k \cdot N \subset I\cdot N,
$$
and, thus,
$$
\al^k\in \ga_k^{-1} \cdot I\cdot N.
$$

We now consider the action of the isometries in the above equation
on the building $X_p$: the group $\Ga$ fixes a point $o\in X_p$,
the images $I\cdot o$ form a finite set. Therefore, the set
$$
\ga_k^{-1} \cdot I\cdot N \cdot o = \ga_k^{-1} \cdot I \cdot o
$$
is bounded in $X$. However, by our assumption, the orbit $\al^k
\cdot o, k\in \N$ is unbounded. Contradiction. \qed

\begin{cor}
Suppose that $\Ga$ has trivial center. Then, if $N<\La$ is a
normal subgroup of $\Ga$ so that $\varphi(N)\subset N$, then $N$
is trivial.
\end{cor}
\proof Since $\Ga$ has trivial center, its only finite normal
subgroups are trivial. \qed

\medskip
We now consider virtual endomorphisms of non-arithmetic lattices.

\begin{prop}
Let $G$ be a rank 1 semisimple Lie group (with finitely many
components), which is not locally isomorphic to $SL(2,\R)$ and
$\Ga<G$ be a non-arithmetic lattice with trivial center. Then for
every virtual endomorphism $\varphi$ of $\Ga$, there exists an
infinite normal subgroup $N<\Ga$ which is $\varphi$-invariant.
\end{prop}
\proof 1. Suppose that $\varphi: \La\to \Ga$ is not injective. Let
$K$ denote the kernel of $\varphi$. Note that this subgroup is not
necessarily normal in $\Ga$. Our assumption that $\Ga$ has trivial
center implies that $\Ga$ acts faithfully on the symmetric space
$X$ associated with $G$ and, hence, the group $K$ is necessarily
infinite.

We first consider the case when $\La$ is normal in $\Ga$. Let
$\ga_1,...,\ga_n$ be the generators of $\Ga$. Consider the
conjugates
$$
K_i:= \ga_i K \ga_i^{-1}\subset \Ga, \quad K_0:= K.
$$
Then
$$
K':=\bigcap_{i=0}^n K_i
$$
is a normal subgroup in $\Ga$. This subgroup is the kernel of the
homomorphism
$$
\Phi: \La\to \prod_{i=0}^n \Ga,
$$
$$
\Phi=(\varphi_0,...,\varphi_n),
$$
where
$$
\varphi_i: \La \to \Ga
$$
is given by
$$
\varphi_i(g)=   \varphi( \ga_i g \ga_i^{-1}), i=1,...,n; \quad
\varphi_0:= \varphi.
$$
We set let $\Ga_i$ denote the $i$-th factor of the product group
$\prod_{i=0}^n \Ga$.

If $K$ is infinite, it contains a free nonabelian subgroup $H$.
Thus, each group $K_i$ contains a free nonabelian subgroup $H_i$,
$i=0,...,n$. Assume, for a moment, that $K'$ is finite.

The group $\La/K'$ embeds in the product group $\prod_{i=0}^n
\Ga$. Therefore, the intersections
$$
\Phi(\La) \cap \Ga_i
$$
contain free nonabelian subgroups (isomorphic to $H_i$),
$i=0,...,n$. Hence, the group $\La/K'$ contains a direct product
of free nonabelian subgroups. This is impossible since the group
$\La/K'$ is isomorphic to a discrete subgroup of isometries of the
negatively curved symmetric space $X$. Contradiction. Thus, $K'$
is infinite and we obtain an infinite normal subgroup $N=K'<\La$
of the group $\Ga$, so that $\varphi(N)=1\subset N$.

\medskip
We now consider the case when $\La$ is not necessarily normal in
$\Ga$. If $\varphi: \La\to \Ga$ is a virtual endomorphism with
infinite kernel, we find a finite index subgroup $\La'\subset \La$
which is normal in $\Ga$. Then the restriction
$\varphi'=\varphi|\La'$ still has infinite kernel and we obtain a
contradiction as above.

\medskip
2. Suppose now that $\varphi$ is injective. Then, by Mostow
rigidity theorem, the homomorphism $\varphi$ is induced by
conjugation via some $\al\in Comm(\Ga)$. Recall that, by the
Margulis theorem, $Comm(\Ga)$ is discrete. Therefore, since $\Ga$,
is a lattice, $\hat\Ga:=Comm(\Ga)$ is a finite extension of $\Ga$.
Therefore, $\La$ has finite index in $\hat\Ga$ and, hence,
contains a finite index subgroup $N<\La$ which is normal in
$\hat\Ga$. In particular,
$$
\al N \al^{-1}=N
$$
and $N<\Ga$ is a normal subgroup. Clearly, $N$ is an infinite
group as required by the proposition. \qed

\medskip
If $G$ is locally isomorphic to $SL(2,\R)$ then the abstract
commensurator of a lattice $\Ga<G$ is not isomorphic to the
commensurator of $\Ga$ in $G$. However, $\Ga$ will contain a
finite index torsion-free subgroup $\La$ which is isomorphic to an
arithmetic subgroup of $G$. Therefore, $\Ga$ will have a virtual
endomorphism $\varphi$ without nontrivial invariant normal
subgroups.

\medskip
By combining the above results we obtain Theorem \ref{main}. \qed

\medskip
We observe that our results should, in principle, generalize to
Gromov-hyperbolic groups which are not lattices. The problem,
however, is that:

1. Among hyperbolic groups $\Ga$, (if we ignore torsion) only
Poincar\'e duality groups (fundamental groups of closed aspherical
manifolds) are known to be {\em weakly cohopfian}, in the sense
that if $\La<\Ga$ is a finite index subgroup and $\varphi: \La\to
\Ga$ is an injective homomorphism, then the image of $\varphi$ is
a finite index subgroup of $\Ga$. Hyperbolic groups which act
geometrically on rank 2 hyperbolic buildings provide good
candidates for weakly cohopfian groups \cite{Bourdon-Pajot(2002)}.
On the other hand, there are no  known examples of 1-ended
hyperbolic groups which are not weakly cohopfian.

\medskip
2. Among hyperbolic groups $\Ga$ which are Poincar\'e duality
groups, the only known examples where the abstract commensurator
$Comm(\Ga)$ is a finite extension of $\Ga$, are the non-arithmetic
lattices. There are few more classes of hyperbolic groups with
small abstract commensurators: Fundamental groups of compact
hyperbolic $n$-manifolds with totally-geodesic boundary, $n\ge 3$,
surface-by-free groups \cite{Farb-Mosher(2002)}, as well as some
rigid examples constructed in \cite{Kapovich-Kleiner00}. (In all
these examples, the whole quasi-isometry group of $\Ga$ is a
finite extension of $\Ga$.) Conjecturally, the fundamental groups
of Gromov-Thurston manifolds \cite{Gromov-Thurston(1987)} should
also have small abstract commensurators.

\begin{question}
Does the action of an arithmetic group $\Ga$ on a rooted tree
constructed in Theorem \ref{main} ever correspond to a
finite-state automaton? Zoran Sunic computed one example (an index
3 subgroup in $PSL(2,\Z)$) where he proved that the number of
states is infinite. See however the examples constructed by
Glasner and Mozes in \cite{Glasner-Mozes(2005)}.
\end{question}

\medskip
{\bf Acknowledgements.} This paper is the result of conversations
the author had with Lucas Sabalka and Zoran Sunic. The author is
indebted to Zoran Sunic for explaining him the basic definitions
of self-similar group actions and providing him with references. I
am also grateful to Benson Farb, Mahan Mj and David Fisher for
remarks and corrections.
%The author was supported by the NSF grant DMS-05-54349.

\bibliography{lit}%C:\texmf\bibtex/bib/lit.bib

\begin{thebibliography}{1}

\bibitem{Bourdon-Pajot(2002)}
{\sc M.~Bourdon and H.~Pajot}, {\em Quasi-conformal geometry and hyperbolic
  geometry}, in Rigidity in dynamics and geometry (Cambridge, 2000), Springer,
  Berlin, 2002, pp.~1--17.

\bibitem{Farb-Mosher(2002)}
{\sc B.~Farb and L.~Mosher}, {\em The geometry of surface-by-free groups},
  Geom. Funct. Anal., 12 (2002), pp.~915--963.

\bibitem{Glasner-Mozes(2005)}
{\sc Y.~Glasner and S.~Mozes}, {\em Automata and square complexes}, Geom.
  Dedicata, 111 (2005), pp.~43--64.

\bibitem{Gromov-Thurston(1987)}
{\sc M.~Gromov and W.~Thurston}, {\em Pinching constants for hyperbolic
  manifolds}, Inventiones Mathematicae, 89 (1987), pp.~1--12.

\bibitem{Kapovich-Kleiner00}
{\sc M.~Kapovich and B.~Kleiner}, {\em Hyperbolic groups with low-dimensional
  boundary}, Ann. Sci. \'Ecole Norm. Sup. (4), 33 (2000), pp.~647--669.

\bibitem{Margulis}
{\sc G.~A. Margulis}, {\em Discrete subgroups of semisimple {L}ie groups},
  Springer-Verlag, Berlin, 1991.

\bibitem{Nekrashevych}
{\sc V.~Nekrashevych}, {\em Virtual endomorphisms of groups}, Algebra Discrete
  Math., 1 (2002), pp.~88--128.

\bibitem{Nekrashevych-Sidky}
{\sc V.~Nekrashevych and S.~Sidky}, {\em Automorphisms of the binary tree:
  {S}tate-closed subgroups and dynamics of 1/2-endomorphisms}, in Groups:
  Topological, Combinatorial and Arithmetic Aspects, T.~W. M\"uller, ed.,
  vol.~311 of LMS Lecture Notes, London Math. Society, 2004, pp.~375--404.

\end{thebibliography}
\bibliographystyle{siam}

\end{document}